\documentclass[12pt]{article}
\usepackage{amsfonts}

\def\hybrid{\topmargin 0pt      \oddsidemargin 0pt
        \headheight 0pt \headsep 0pt
        \voffset=-0.5cm
        \textwidth 6.5in        
        \textheight 9in         
        \marginparwidth 0.0in
        \parskip 5pt plus 1pt   \jot = 1.5ex}
\catcode`\@=11
\def\marginnote#1{}

\newcount\hour
\newcount\minute
\newtoks\amorpm
\hour=\time\divide\hour by60
\minute=\time{\multiply\hour by60 \global\advance\minute by-\hour}
\edef\standardtime{{\ifnum\hour<12 \global\amorpm={am}%
        \else\global\amorpm={pm}\advance\hour by-12 \fi
        \ifnum\hour=0 \hour=12 \fi
        \number\hour:\ifnum\minute<10 0\fi\number\minute\the\amorpm}}
\edef\militarytime{\number\hour:\ifnum\minute<10 0\fi\number\minute}

\def\draftlabel#1{{\@bsphack\if@filesw {\let\thepage\relax
   \xdef\@gtempa{\write\@auxout{\string
      \newlabel{#1}{{\@currentlabel}{\thepage}}}}}\@gtempa
   \if@nobreak \ifvmode\nobreak\fi\fi\fi\@esphack}
        \gdef\@eqnlabel{#1}}
\def\@eqnlabel{}
\def\@vacuum{}
\def\draftmarginnote#1{\marginpar{\raggedright\scriptsize\tt#1}}
\def\draftlabel#1{{\@bsphack\if@filesw {\let\thepage\relax
   \xdef\@gtempa{\write\@auxout{\string
      \newlabel{#1}{{\@currentlabel}{\thepage}}}}}\@gtempa
   \if@nobreak \ifvmode\nobreak\fi\fi\fi\@esphack}
        \gdef\@eqnlabel{#1}}
\def\@eqnlabel{}
\def\@vacuum{}
\def\draftmarginnote#1{\marginpar{\raggedright\scriptsize\tt#1}}

\def\draft{\oddsidemargin -.5truein
        \def\@oddfoot{\sl preliminary draft \hfil
        \rm\thepage\hfil\sl\today\quad\militarytime}
        \let\@evenfoot\@oddfoot \overfullrule 3pt
        \let\label=\draftlabel
        \let\marginnote=\draftmarginnote
   \def\@eqnnum{(\theequation)\rlap{\kern\marginparsep\tt\@eqnlabel}%
\global\let\@eqnlabel\@vacuum}  }

\def\numberbysection{\@addtoreset{equation}{section}
        \def\theequation{\thesection.\arabic{equation}}}

\def\underline#1{\relax\ifmmode\@@underline#1\else
        $\@@underline{\hbox{#1}}$\relax\fi}

\def\titlepage{\@restonecolfalse\if@twocolumn\@restonecoltrue\onecolumn
     \else \newpage \fi \thispagestyle{empty}\c@page\z@
        \def\thefootnote{\fnsymbol{footnote}} }

\def\endtitlepage{\if@restonecol\twocolumn \else  \fi
        \def\thefootnote{\arabic{footnote}}
        \setcounter{footnote}{0}}  
\relax

\numberbysection
\hybrid

\def\beq{\begin{equation}}
\def\eeq{\end{equation}}
\def\bea{\begin{eqnarray}}
\def\eea{\end{eqnarray}}
\def\p{\partial}
\def\G{\Gamma}
\def\g{\gamma}

\def\L{{\cal L}}

\def\a{\alpha}
\def\b{\beta}
\def\e{\varepsilon}
\def\l{\lambda}
\def\f{\varphi}

\def\A{{\cal A}}

\def\V{{\cal V}}
\def\D{{\cal D}}
\def\F{{\cal F}}

\def\L{{\cal L}}

\def\O{{\cal O}}

\def\dim{{\rm dim}}
\def\res{{\rm res}}

\def \matrix #1 {\left(\begin{array}{cc} #1 \end{array}\right)}
\newtheorem{theo}{Theorem}[section]

\newtheorem{cor}{Corollary}[section]
\newtheorem{lem}{Lemma}[section]

\begin{document}
\begin{titlepage}
\title{Integrable linear equations and the Riemann-Schottky problem}

\author{I.Krichever \thanks{Columbia University, New York, USA and
Landau Institute for Theoretical Physics, Moscow, Russia; e-mail:
krichev@math.columbia.edu. Research is supported in part by National Science
Foundation under the grant DMS-04-05519.}}

\date{April 8, 2005\,\footnote {the revised version November 30,2005} }

\maketitle

\begin{abstract} We prove that an indecomposable principally polarized
abelian variety $X$ is the Jacobain of a curve if and only if there exist vectors
$U\neq 0,V$ such that the roots $x_i(y)$ of the theta-functional equation
$\theta(Ux+Vy+Z)=0$
satisfy the equations of motion of the {\it formal infinite-dimensional
Calogero-Moser system}.
\end{abstract}

\end{titlepage}

\section{Introduction}
The Riemann-Schottky problem on the characterization of the Jacobians of curves among
abelian
varieties is more that 120 years old. Quite a few geometrical characterizations
of the Jacobians have been found. None of them provides an explicit system of equations for the image of the Jacobian locus
in the projective space under the level two theta imbedding.

The first effective solution of the Riemann-Schottky problem was obtained by
T.Shiota (\cite{shiota}), who proved the famous Novikov's conjecture:

{\it An indecomposable principally polarized abelian variety $(X,\theta)$
is the Jacobian of a curve of a genus g if and only if there exist $g$-dimensional vectors
$U\neq 0,V,W$ such that the function
\beq\label{u}
u(x,y,t)=-2\p_x^2 \ln \theta (Ux+Vy+Wt+Z)
\eeq
is a solution of the Kadomtsev-Petviashvilii (KP) equation}
\beq\label{kp}
{3}u_{yy}=\left(4u_t+6uu_x-u_{xxx}\right)_x .
\eeq
Here $\theta(Z)=\theta(Z|B)$ is the Riemann theta-function,
\beq\label{teta1}
\theta(z)=\sum_{m\in \mathbb{Z}^g} e^{2\pi i(z,m)+\pi i(Bm,m)},\ \ (z,m)=m_1z_1+\cdots+m_gz_g,
\eeq
where $B$ is the corresponding symmetric matrix with a positive definite imaginary part.

It is easy to show (\cite{dub}) that the KP equation with $u$ of the form (\ref{u})
is in fact equivalent to the following system of algebraic equations for the fourth order
derivatives of the level two theta constants:
\beq\label{dub}
\p_U^4\Theta[\e,0]-\p_U\p_W\Theta[\e,0]+\p_V^2\Theta[\e,0]+c\Theta[\e,0]=0, \ c={\rm const}.
\eeq
Here $\Theta[\e,0]=\Theta[\e,0](0)$, where $\Theta[\e,0](z)=\theta[\e,0](2z|2B)$
are level two theta-functions with half-integer characteristics $\e\in {1\over 2} \mathbb{Z}_2^g$.

The KP equation admits the so-called zero-curvature representation (\cite{zakh,dr})
which is the compatibility condition for the following over-determined
system of linear equations:
\bea
{}&\left(\p_y-\p_x^2+u\right)\psi=0\,\label{lax1},\\
{}&\left(\p_t-\p_x^3+{3\over 2}\,\p_x+w\right)\psi=0\,.\label{lax2}
\eea
The main goal of the present paper is to show that the KP equation contains an excessive
information and that the Jacobians can be characterized in terms of {\it only the first}
of its auxiliary linear equations.

\begin{theo} An indecomposable principally polarized abelian variety $(X,\theta)$
is the Jacobian of a curve of genus g if and only if there exist $g$-dimensional vectors
$U\neq 0, V,A $ such that equation (\ref{lax1}) is satisfied for
\beq\label{u2}
u=-2\p_x^2 \ln \theta (Ux+Vy+Z)
\eeq
and
\beq\label{p}
\psi={\theta(A+Ux+Vy+Z)\over \theta(Ux+Vy+Z)}\, e^{p\,x+E\,y},
\eeq
where $p,E$ are constants.
\end{theo}
The ``if\," part of this statement follows from the exact theta-functional expression
for the Baker-Akhiezer function (\cite{kr1,kr2}).

The addition formula for the Riemann theta-function directly implies that
equation (\ref{lax1}) with $u$ and $\psi$ of the form (\ref{u2}) and (\ref{p}) is
equivalent to the system of equations.
\beq\label{gr}
\left(\p_V-\p_U^2-2p\,\p_U+(E-p^2)\right)\, \Theta[\e,0](A/2)=0,\ \ \e\in {1\over 2}
\mathbb{Z}_2^g.
\eeq
Recently Theorem 1.1 was proved by E.Arbarello, G. Marini and the author (\cite{flex})
under the additional assumption that the closure $\langle A\rangle$ of the subgroup  of $X$
generated by $A$ is irreducible. The geometric interpretation of Theorem 1.1
is equivalent to the characterization of the Jacobians via flexes of Kummer varieties
(see details in \cite{flex}), which is a particular case of
the so-called {\it trisecant conjecture},
first formulated in \cite{wel1}.

Theorem 1.1 is not the strongest form of our main result. What we really prove is that
the Jacobian locus in the space of principally polarized abelian varieties is characterized
by a system of equations which formally can be seen as the equations of motion
of  the {\it infinite-dimensional} Calogero-Moser system.

Let $\tau(x,y)$ be an entire function of the complex variable $x$ smoothly depending on
a parameter $y$. Consider the equation
\beq\label{rr}
\res_x\left(\p_y^2\ln \tau+2\left(\p_x^2\ln \tau\right)^2\right)=0,
\eeq
which means that the meromorphic function given by the left hand side of (\ref{rr})
has no {\it residues} in the $x$ variable. If $x_i(y)$ is a simple zero of $\tau$, i.e. $\tau(x_i(y),y)=0,
\p_x\tau(x_i(y),y)\neq 0$, then (\ref{rr}) implies
\beq\label{cm5}
\ddot x_i=2w_i,\
\eeq
where ``dots" stands for the $y$-derivatives and $w_i$ is the third coefficient of the Laurent
expansion of $u(x,y)=-2\p_x^2 \tau(x,y)$ at $x_i$, i.e.
\beq\label{e1}
u(x,y)={2\over (x-x_i(y))^2}+v_i(y)+w_i(y)(x-x_i(y))+\cdots
\eeq
Formally, if we represent $\tau$ as an infinite product,
\beq\label{roots}
\tau(x,y)=c(y)\prod_i(x-x_i(y)),
\eeq
then equation (\ref{rr}) can be written as the infinite system of equations
\beq\label{cm6}
\ddot x_i=-4\sum_{j\neq i} {1\over (x_i-x_j)^3}\,.
\eeq
Equations (\ref{cm6}) are purely formal because, even if $\tau$ has simple zeros at $y=0$,
then in the general case there is no nontrivial interval in $y$ where the zeros stay simple.
At the moment, the only reason for representing (\ref{cm5}) in the form
(\ref{cm6}) is to show that in the case when $\tau$ is a rational, trigonometric or
elliptic polynomial  the system (\ref{cm5}) coincides with the equations of
motion for the rational, trigonometrical or elliptic Calogero-Moser systems, respectively.

Equations (\ref{cm5}) for the zeros of the function
$\tau=\theta(Ux+Vy+Z)$ were derived in \cite{flex} as a direct corollary of the
assumptions of Theorem 1.1.
Simple expansion of $\theta$ at the points of its divisor $z\in\Theta: \theta(z)=0$
gives the equation
\beq\label{cm7}
[(\p_2\theta)^2-(\p_1^2\theta)^2]\p_1^2\theta
+2[\p_1^2\theta\p_1^3\theta-\p_2\theta\p_1\p_2\theta]\p_1\theta+
[\p_2^2\theta-\p_1^4\theta](\p_1\theta)^2=0\ ({\rm mod}\, \theta)
\eeq
{\it which  is valid on $\Theta$.} Here and below $\Theta$ is the divisor on $X$ defined by
the equation $\theta(Z)=0$ and $\p_1$ and $\p_2$ are constant vector field
on $\mathbb{C}^g$ corresponding to the vectors $U$ and $V$.

It would be very interesting to understand if any reasonable
general theory of equations (\ref{rr}) exists. The following form of our main
result shows that in any case such theory has to be interesting and non-trivial.

Let $\Theta_1$ be defined by the equations
$\Theta_1=\{Z:\theta(Z)=\partial_1\theta(Z)=0\}$. The $\p_1$-invariant subset $\Sigma$
of $\Theta_1$ will be called the {\it singular locus}.

\begin{theo} An indecomposable principally polarized abelian variety $(X,\theta)$
is the Jacobian of a curve of genus g if and only if there exist $g$-dimensional vectors
$U\neq 0,V$, such that for each $Z\in \mathbb C^g\setminus \Sigma$
equation (\ref{rr}) for the function $\tau(x,y)=\theta (Ux+Vy+Z)$ is satisfied, i.e.
equation (\ref{cm7}) is valid on $\Theta$.
\end{theo}
The main idea of Shiota's proof of the Novikov's conjecture is to
show that if $u$ is as in (\ref{u}) and satisfies the KP equation, then
it can be extended to a $\tau$-function of the KP {\it hierarchy},
as a {\it global} holomorphic function of the infinite number of variables
$t=\{t_i\},\ t_1=x, t_2=y, t_3=t$. Local existence of $\tau$ directly follows from
the KP equation. The global existence of the $\tau$-function
is crucial. The rest is a corollary of the KP theory and the theory of commuting
ordinary differential operators developed by Buchnall-Chaundy (\cite{ch1,ch2}) and
the author (\cite{kr1,kr2}).

The core of the problem is that there is a homological obstruction for the global existence
of $\tau$. It is controlled by the cohomology group $H^1(\mathbb C^g\setminus \Sigma, \V)$,
where $\V$ is the sheaf of $\partial_1$-{\it invariant} meromorphic
functions on $\mathbb{C}^g\setminus \Sigma$ with poles along $\Theta$ (see details in \cite{ac}).
The hardest part of the Shiota's work (clarified in \cite{ac}) is the proof that
the locus $\Sigma$ is empty. That insures vanishing of $H^1(\mathbb C^g, \V)$.
Analogous obstructions have occurred in all the other attempts to apply the theory
of soliton equations
to various characterization problems in the theory of abelian varieties. None of
them had been completely successful. Only partial results were obtained.
(Note that Theorem 1.1 in one of its equivalent forms was proved earlier in \cite{mar}
under the additional assumption that $\Theta_1$ does not contain $\p_1$-invariant
line.)

Strictly speaking, the KP equation and the KP hierarchy are not used in the present paper.
But our main construction of {\it the formal wave solutions} of (\ref{lax1})
is reminiscent to the construction of the $\tau$-function. All its difficulties
can be traced back to those in the Shiota's work.
The wave solution of (\ref{lax1}) is a solution of the form
\beq\label{ps}
\psi(x,y,k)=e^{kx+(k^2+b)y}\left(1+\sum_{s=1}^{\infty}\xi_s(x,y)\,k^{-s}\right)\,.
\eeq
At the beginning of the next section we show that the assumptions of Theorem 1.2 are
necessary and sufficient conditions for the {\it local} existence of the wave solutions
such that
\beq\label{xi22}
\xi_s={\tau_s(Ux+Vy+Z,y)\over\theta (Ux+Vy+Z)}\ , \ Z\notin \Sigma,
\eeq
where $\tau_s(Z,y)$, as a function of $Z$, is holomorphic in some open domain in $\mathbb C^g$.
The functions $\xi_s$ are defined recurrently by the equation $2\partial_1\xi_{s+1}=
\p_y \xi_s-\partial_1^2\xi_s+u\xi_s$. Therefore, the global existence of $\xi_s$ is controlled
by the same cohomology group $H^1(\mathbb C\setminus \Sigma, \V)$ as above. At the local level
the main problem is to find translational invariant normalization of $\xi_s$ which defines
wave solutions uniquely up  to a $\p_1$-invariant factor.

In the case of periodic potentials $u(x+T,y)=u(x)$ the normalization problem for the
wave functions was solved D.Phong and the author in \cite{kp2}.
It was shown that the condition that $\xi_s$ {\it is periodic} completely determines
the $y$-dependence of the integration constants and the corresponding wave solutions are
related by a $x$-independent factor. In general, the potential $u=-2\p_x^2\theta(Ux+Vy+Z)$ is
only quasi-periodic in $x$. In that case the solution of the normalization problem is technically
more involved but mainly goes along the same lines as in the periodic case.
The corresponding wave solutions are called $\l$-periodic.

In the last section we show that for each $Z\notin \Sigma$ a local
$\l$-periodic wave solution is the common eigenfunction of a commutative ring $\A^Z$ of
ordinary differential operators. The coefficients of these operators are independent of
ambiguities in the construction of $\psi$. For generic $Z$ the ring $A^Z$ is maximal
and the corresponding spectral curve $\G$ is $Z$-independent. The correspondence
$j:Z\longmapsto \A^Z$ allows us to make the next
crucial step and prove the global existence of the wave function.
Namely, on $(X\setminus\Sigma)$ the wave function can be globally defined as
the preimage $j^*\psi_{BA}$ under $j$ of the Baker-Akhiezer function on $\G$ and
then can be extended on $X$ by usual Hartorg's arguments. The global existence of the
wave function implies that $X$ contains an orbit of the
KP hierarchy, as an abelian subvariety. The orbit is isomorphic to the generalized
Jacobian $J(\G)={\rm Pic}^0(\G)$ of the spectral curve (\cite{shiota}).
Therefore, the generalized Jacobian
is compact. The compactness of ${\rm Pic}^0(\G)$ implies that the spectral curve is smooth
and the correspondence $j$ extends by linearity and define an isomorphism
$j: X\to J(\G)$.

\section{$\l$-periodic wave solutions}
As it was mentioned above, formal Calogero-Moser equations (\ref{cm5}) were derived
in \cite{flex} as a necessary condition for the existence of a meromorphic solution
to equation (\ref{lax1}).

Let $\tau(x,y)$ be a holomorphic function of the variable $x$ in some open domain
$D\in \mathbb C$ smoothly depending on a parameter $y$. Suppose that for each $y$
the zeros of $\tau$ are simple,
\beq\label{xi}
\tau(x_i(y),y)=0, \tau_x(x_i(y),y)\neq 0.
\eeq
\begin{lem} (\cite{flex})
If equation (\ref{lax1}) with the potential
$u=-2\p_x^2\ln \tau(x,y)$ has a meromorphic in $D$ solution $\psi_0(x,y)$, then
equations (\ref{cm5}) hold.
\end{lem}
{\it Proof.}
Consider the Laurent expansions of $\psi_0$ and $u$ in the neighborhood of one of the zeros
$x_i$ of $\tau$:
\beq\label{ue}
u={2\over (x-x_i)^2}+v_i+w_i(x-x_i)+\ldots
\eeq
\beq\label{psie}
\psi_0={\a_i\over x-x_i}+\b_i+\g_i(x-x_i)+\delta_i(x-x_i)^2+\ldots
\eeq
(All coefficients in these expansions are smooth functions of the variable $y$).
Substitution of (\ref{ue},\ref{psie}) in (\ref{lax1})
gives a system of equations. The first three of them
are
\beq\label{eq1}
\a_i \dot x_i+2\b_i=0,
\eeq
\beq\label{eq2}
\dot\a_i+\a_i v_i+2\g_i=0,
\eeq
\beq\label{eq3}
\dot\b_i+v_i\b_i-\g_i\dot x_i+\a_i w_i=0.
\eeq
Taking the $y$-derivative of the first equation and using two others we get
(\ref{cm5}).

Let us show that equations (\ref{cm5}) are sufficient for the existence of
meromorphic wave solutions.
\begin{lem}
Suppose that equations (\ref{cm5}) for the zeros of $\tau(x,y)$ hold. Then there exist
meromorphic wave solutions of equation (\ref{lax1}) that have simple poles at $x_i$ and
are holomorphic everywhere else.
\end{lem}
{\it Proof.} Substitution of (\ref{ps}) into (\ref{lax1}) gives a recurrent system of
equations
\beq\label{xis}
2\xi_{s+1}'=\p_y\xi_s+u\xi_s-\xi_s''
\eeq
We are going to prove by induction that this system has meromorphic solutions with
simple poles at all the zeros $x_i$ of $\tau$.

Let us expand $\xi_s$ at $x_i$:
\beq\label{5}
\xi_s={r_s\over x-x_i}+r_{s0}+r_{s1}(x-x_i)\,,
\eeq
where for brevity we omit the index $i$ in the notations for the coefficients of this
expansion.
Suppose that $\xi_s$ are defined and equation (\ref{xis}) has a meromorphic solution.
Then the right hand side of (\ref{xis}) has the zero residue at $x=x_i$, i.e.,
\beq\label{res}
{\rm res}_{x_i}\left(\p_y\xi_s+u\xi_s-\xi_s''\right)=\dot r_s+v_ir_s+2r_{s1}=0
\eeq
We need to show that the residue of the next equation vanishes also.
From (\ref{xis}) it follows that the coefficients of the Laurent expansion for $\xi_{s+1}$
are equal to
\beq\label{6}
r_{s+1}=-\dot x_ir_s-2r_{s0},
\eeq
\beq\label{7}
 2r_{s+1,1}=\dot r_{s0}-r_{s1}+w_ir_s+v_ir_{s0}\,.
\eeq
These equations imply
\beq
\dot r_{s+1}+v_ir_{s+1}+2r_{s+1,1}=-r_s(\ddot x_i-2w_i)-\dot x_i(\dot r_s-v_ir_ss+2r_{s1})=0,
\eeq
and the lemma is proved.

Our next goal is to fix a {\it translation-invariant} normalization of $\xi_s$
which defines wave functions uniquely up to a $x$-independent factor.
It is instructive to consider first the case of the periodic potentials $u(x+1,y)=u(x,y)$
(see details in \cite{kp2}).

Equations (\ref{xis}) are solved recursively by the formulae
\beq
\xi_{s+1}(x,y)=c_{s+1}(y)+\xi_{s+1}^0(x,y)\,,\label{kp1}
\eeq
\beq\label{kp2}
\xi_{s+1}^0(x,y)={1\over 2}\int_{x_0}^x (\p_y \xi_s-\xi_s''+u\xi_s)\,dx\, ,
\eeq
where $c_s(y)$ are {\it arbitrary} functions of the variable $y$.
Let us show that the periodicity condition $\xi_s(x+1,y)=\xi_s(x,y)$
defines the functions $c_s(y)$ uniquely up to an additive constant.
Assume that $\xi_{s-1}$ is known  and satisfies the condition that the corresponding
function $\xi_s^0$ is periodic.
The choice of the function $c_s(y)$ does not affect the periodicity property of
$\xi_s$, but it does affect the periodicity in $x$ of the function
$\xi_{s+1}^0(x,y)$. In order to make  $\xi_{s+1}^0(x,y)$ periodic,
the function $c_s(y)$ should satisfy the linear differential equation
\beq\label{kp4}
\p_y c_s(y)+B(y)\,c_s(y)+\int_{x_0}^{x_0+1} \left(\p_y \xi_s^0(x,y)+
u(x,y)\,\xi_s^0(x,y)\right)\,dx\ ,
\eeq
where $B(y)=\int_{x_0}^{x_0+1} u\, dx$.
This defines $c_s$ uniquely up to a constant.

In the general case, when $u$ is quasi-periodic, the normalization of the wave functions
is defined along the same lines.

Let $Y_U=\langle Ux\rangle$ be the closure of the group
$Ux$ in $X$. Shifting $Y_U$ if needed, we may assume, without loss of generality, that
$Y_U$ is not in the singular locus, $Y_U\notin\Sigma$. Then, for a sufficiently small
$y$, we have $Y_U+Vy\notin\Sigma$ as well.
Consider the restriction of the theta-function onto the affine subspace
$\mathbb C^d+Vy$, where $\mathbb C^d=\pi^{-1}(Y_U)$, and
$\pi: \mathbb C^g\to X=\mathbb C^g/\Lambda$ is the universal cover of $X$:
\beq\label{ttt1}
\tau (z,y)=\theta(z+Vy), \ \ z\in \mathbb C^d.
\eeq
The function $u(z,y)=-2\partial_1^2\ln \tau$ is periodic with respect to the lattice
$\Lambda_U=\Lambda\cap \mathbb C^d$ and, for fixed $y$, has a double pole along the divisor
$\Theta^{\,U}(y)=\left(\Theta-Vy\right)\cap \mathbb C^d$.

\begin{lem} Let equation (\ref{rr}) for $\tau(Ux+z,y)$ hold and let
$\l$ be a vector of the sublattice $\Lambda_U=\Lambda\cap \mathbb C^d\subset \mathbb C^g$. Then:

(i) equation (\ref{lax1}) with the potential $u(Ux+z,y)$
has a wave solution of  the form $\psi=e^{kx+k^2y}\phi(Ux+z,y,k)$
such that the coefficients $\xi_s(z,y)$ of the formal series
\beq\label{psi2}
\phi(z,y,k)=e^{by}\left(1+\sum_{s=1}^{\infty}\xi_s(z,y)\, k^{-s}\right)
\eeq
are $\l$-periodic meromorphic functions of the variable $z\in \mathbb C^d$
with a simple pole at
the divisor $\Theta^U(y)$,
\beq\label{v1}
\xi_s(z+\l,y)=\xi_s(z,y)={\tau_s(z,y)\over \tau(z,y)}\, ;
\eeq

(ii) $\phi(z,y,k)$ is unique up to a factor $\rho(z,k)$ that is $\partial_1$-invariant
and holomorphic in $z$,
\beq\label{v2}
\phi_1(z,y,k)=\phi(z,y,k)\rho(z,k), \ \partial_1\rho=0.
\eeq
\end{lem}
{\it Proof.} The functions $\xi_s(z)$ are defined recursively by the equations
\beq\label{xis1}
2\partial_1\xi_{s+1}=\p_y\xi_s+(u+b)\xi_s-\partial_1^2\xi_s.
\eeq
A particular solution of the first equation $2\partial_1\xi_1=u+b$ is given by the formula
\beq\label{v5}
2\xi_1^0=-2\partial_1\ln \tau +(l,z)\ b,
\eeq
where $(l,z)$ is a linear form on $\mathbb C^d$ given by the scalar product of $z$ with
a vector $l\in \mathbb C^d$ such that $(l,U)=1$, and $(l,\l)=1.$ The periodicity
condition for $\xi_1^0$ defines the constant $b$
\beq\label{v6}
b=2\partial_1\ln \tau(z+\l,y)-2\partial_1\ln \tau(z,y)\,,
\eeq
which depends only on a choice of the lattice vector $\l$.
A change of the potential by an additive constant does not affect the results of the
previous lemma. Therefore, equations (\ref{cm5}) are sufficient for the local solvability
of (\ref{xis1}) in any domain, where $\tau(z+Ux,y)$ has simple zeros, i.e., outside
of the set $\Theta_1^{\,U}(y)=\left(\Theta_1-Vy\right)\cap \mathbb C^d$.
Recall that $\Theta_1=\Theta\cap \partial_1\Theta$.
This set  does not contain a $\partial_1$-invariant line because any such line is dense in
$Y_U$. Therefore, the sheaf $\V_0$ of $\partial_1$-invariant meromorphic functions on
$\mathbb{C}^d\setminus \Theta^{\,U}_1(y)$ with poles along the divisor $\Theta^{\,U}(y)$
coincides with the sheaf of holomorphic $\partial_1$-invariant functions.
That implies the vanishing of $H^1({C}^d\setminus \Theta^{\,U}_1(y), \V_0)$ and the
existence of global meromorphic solutions $\xi_s^0$ of (\ref{xis1}) which have a simple
pole at the divisor $\Theta^{\,U}(y)$ (see details in \cite{shiota, ac}).
If $\xi_s^0$ are fixed, then the general global meromorphic solutions are given by the
formula $\xi_s=\xi_s^0+c_s$, where the constant of integration $c_s(z,y)$
is a holomorphic $\partial_1$-invariant function of the variable $z$.

Let us assume, as in the example above, that a $\l$-periodic solution
$\xi_{s-1}$ is known  and that it satisfies the condition that there exists a periodic
solution $\xi_s^0$ of the next equation. Let $\xi_{s+1}^*$ be a solution of
(\ref{xis1}) for fixed $\xi_s^0$. Then it is easy to see that the function
\beq\label{v7}
\xi_{s+1}^0(z,y)=\xi_{s+1}^*(z,y)+c_s(z,y)\,\xi_1^0(z,y)+{(l,z)\over 2}\p_y c_s (z,y),
\eeq
is a solution of (\ref{xis1}) for $\xi_s=\xi_s^0+c_s$.
A choice of a $\l$-periodic $\partial_1$-invariant function $c_s(z,y)$ does not affect the
periodicity property of $\xi_s$, but it does affect the periodicity of the function
$\xi_{s+1}^0$. In order to make  $\xi_{s+1}^0$ periodic,
the function $c_s(z,y)$ should satisfy the linear differential equation
\beq\label{kp41}
\p_y c_s(z,y)=2\xi_{s+1}^*(z+\l,y)-2\xi_{s+1}^*(z,y)\,.
\eeq
This equation, together with an initial condition $c_s(z)=c_s(z,0)$ uniquely defines $c_s(x,y)$.
The induction step is then completed. We have shown that the ratio of two
periodic formal series $\phi_1$ and $\phi$ is $y$-independent. Therefore,
equation (\ref{v2}), where $\rho(z,k)$ is defined by the evaluation of
the both sides at $y=0$, holds. The lemma is thus proven.

\begin{cor} Let $\l_1,\ldots,\l_d$ be a set of linear independent vectors of
the lattice $\Lambda_U$ and let $z_0$ be a point of $\mathbb C^d$. Then, under
the assumptions of the previous lemma,
 there is a unique wave solution of equation (\ref{lax1}) such that
the corresponding formal series $\phi(z,y,k;z_0)$ is quasi-periodic with respect
to $\Lambda_U$, i.e. for
$\l\in \Lambda_U$
\beq\label{v10}
\phi(z+\l,y,k;z_0)=\phi(z,y,k;z_0)\,\mu_{\l}(k)
\eeq
and satisfies the normalization conditions
\beq\label{v11}
\mu_{\l_i}(k)=1,\ \ \ \phi(z_0,0,k;z_0)=1.
\eeq
\end{cor}
The proof is identical to that of the part (b) of the Lemma 12 in \cite{shiota}.
Let us briefly present its main steps. As shown above, there exist wave solutions
corresponding to $\phi$ which are $\l_1$-periodic. Moreover, from the statement (ii) above
it follows that for any $\l'\in \Lambda_U$
\beq\label{v12}
\phi(z+\l,y,k)=\phi(z,y,k)\,\rho_{\l}(z,k)\,,
\eeq
where the coefficients of $\rho_{\l}$ are $\partial_1$-invariant holomorphic functions. Then
the same arguments as in \cite{shiota} show that there exists a $\partial_1$-invariant series
$f(z,k)$ with holomorphic in $z$ coefficients and formal series $\mu_{\l}^0(k)$ with constant
coefficients such that the equation
\beq\label{v13}
f(z+\l,k)\rho_{\l}(z,k)=f(z,k)\,\mu_{\l}(k)
\eeq
holds. The ambiguity in the choice of $f$ and $\mu$ corresponds to the multiplication by
the exponent of a linear form in $z$ vanishing on $U$, i.e.
\beq\label{v14}
f'(z,k)=f(z,k)\,e^{(b(k),\,z)}, \ \ \mu_{\l}'(k)=\mu_{\l}(k)\,e^{(b(k),\,\l)}, \ \ (b(k),\,U)=0,
\eeq
where $b(k)=\sum_sb_s k^{-s}$ is a formal series with vector-coefficients that are
orthogonal to $U$. The vector $U$ is in general position with respect to the lattice. Therefore,
the ambiguity can be uniquely fixed by imposing $(d-1)$ normalizing conditions $\mu_{\l_i}(k)=1,\
i>1$ (recall that $\mu_{\l_1}(k)=1$ by construction).

The formal series $f\phi$ is quasi-periodic and its multiplicators satisfy
(\ref{v11}). Then,  by that properties it is defined uniquely up to a factor
which is constant in $z$ and $y$.
Therefore, for the unique definition of $\phi_0$ it is enough to fix its evaluation at
$z_0$ and $y=0$. The corollary is proved.

\section{The spectral curve}

In this section we show that $\l$-periodic wave solutions of equation (\ref{lax1}), with
$u$ as in (\ref{u2}), are common eigenfunctions of rings of commuting operators
and identify $X$ with the Jacobian of the spectral curve of these rings.

Note that a simple shift $z\to z+Z$, where $Z\notin \Sigma,$ gives
$\l$-periodic wave solutions with meromorphic coefficients along the affine
subspaces $Z+\mathbb C^d$. Theses $\lambda$-periodic wave solutions are related to each other
by $\partial_1$-invariant factor. Therefore choosing, in the neighborhood of any
$Z\notin \Sigma,$ a hyperplane orthogonal to the vector $U$ and
fixing initial data on this hyperplane at $y=0,$ we define the corresponding
series $\phi(z+Z,y,k)$ as a {\it local} meromorphic function of $Z$ and the
{\it global} meromorphic function of $z$.

\begin{lem} Let the assumptions of Theorem 1.2 hold. Then there is a unique
pseudo-differential operator
\beq\label{LL}
\L(Z,\p_x)=\p_x+\sum_{s=1}^{\infty} w_s(Z)\p_x^{-s}
\eeq
such that
\beq\label{kk}
\L(Ux+Vy+Z,\p_x)\,\psi=k\,\psi\,,
\eeq
where $\psi=e^{kx+k^2y} \phi(Ux+Z,y,k)$ is a $\l$-periodic solution of
(\ref{lax1}).
The coefficients $w_s(Z)$ of $\L$  are meromorphic functions on the abelian variety $X$
with poles along  the divisor $\Theta$.
\end{lem}
{\it Proof.} The construction of $\L$ is standard for the KP theory. First we define
$\L$ as a pseudo-differential operator with coefficients $w_s(Z,y)$,
which are functions of $Z$ and $y$.

Let $\psi$ be a
$\l$-periodic wave solution. The substitution of (\ref{psi2}) in (\ref{kk})
gives a system of equations
that recursively define $w_s(Z,y)$ as differential polynomials in $\xi_s(Z,y)$.
The coefficients of $\psi$ are local meromorphic functions of $Z$, but
the coefficients of $\L$ are well-defined
{\it global meromorphic functions} of on $\mathbb C^g\setminus\Sigma$, because
different $\l$-periodic wave solutions are related to each other by $\partial_1$-invariant
factor, which does not affect $\L$. The singular locus is
of codimension $\geq 2$. Then Hartog's holomorphic extension theorem implies that
$w_s(Z,y)$ can be extended to a global meromorpic function on $\mathbb C^g$.

The translational invariance of $u$ implies the translational invariance of
the $\l$-periodic wave solutions. Indeed, for any constant $s$ the series
$\phi(Vs+Z,y-s,k)$ and $\phi(Z,y,k)$ correspond to $\l$-periodic solutions
of the same equation. Therefore, they coincide up to a $\p_1$-invariant factor.
This factor does not affect $\L$. Hence, $w_s(Z,y)=w_s(Vy+Z)$.

The $\l$-periodic wave functions corresponding to $Z$ and
$Z+\lambda'$ for any $\lambda'\in \Lambda$
are also related to each other by a $\partial_1$-invariant factor:
\beq\label{tr1}
\p_1\left(\phi_1(Z+\l',y,k)\phi^{-1}(Z,y,k)\right)=0.
\eeq
Hence, $w_s$ are periodic with respect to $\Lambda$ and therefore are
meromorphic functions on the abelian variety $X$.
The lemma is proved.

Consider now the differential parts of the pseudo-differential operators $\L^m$.
Let $\L^m_+$ be the differential operator such that
$\L^m_-=\L^m-\L^m_+=F_m\p^{-1}+O(\p^{-2})$. The leading
coefficient $F_m$ of $\L^m_-$ is the residue of $\L^m$:
\beq\label{res1}
F_m={\rm res}_{\p}\  \L^m.
\eeq
From the construction of $\L$ it follows that $[\p_y-\p^2_x+u, \L^n]=0$. Hence,
\beq\label{lax}
[\p_y-\p_x^2+u,\L^m_+]=-[\p_y-\p_x^2+u, \L^m_-]=2\p_x F_m.
\eeq
The functions $F_m$ are differential polynomials in the coefficients $w_s$ of $\L$.
Hence, $F_m(Z)$ are meromorphic functions on $X$. Next statement is crucial for
the proof of the existence of commuting differential operators associated with $u$.
\begin{lem} The abelian functions $F_m$ have at most the second order pole on the divisor
$\Theta$.
\end{lem}
{\it Proof.} We need a few more standard constructions from the KP theory.
If $\psi$ is as in Lemma 3.1, then  there exists a unique pseudo-differential
operator $\Phi$ such that
\beq\label{S}
\psi=\Phi e^{kx+k^2y},\ \ \Phi=1+\sum_{s=1}^{\infty}\f_s(Ux+Z,y)\p_x^{-s}.
\eeq
The coefficients of $\Phi$ are universal differential polynomials on $\xi_s$.
Therefore, $\f_s(z+Z,y)$ is a global meromorphic function of $z\in C^d$ and
a local meromorphic function of $Z\notin \Sigma$.
Note that $\L=\Phi(\p_x)\, \Phi^{-1}$.

Consider the dual wave function defined by the left action of the operator $\Phi^{-1}$:
$\psi^+=\left(e^{-kx-k^2y}\right)\Phi^{-1}$.
Recall that the left action of a pseudo-differential operator is the formal adjoint action
under which the left action of $\p_x$ on a function $f$ is $(f\p_x)=-\p_xf$.
If $\psi$ is a formal wave solution of (\ref{lax}),
then $\psi^+$ is a solution of the adjoint equation
\beq\label{adj}
(-\p_y-\p_x^2+u)\psi^+=0.
\eeq
The same arguments, as before, prove that if equations (\ref{cm5}) for poles of $u$
hold then $\xi_s^+$ have simple poles at the poles of $u$. Therefore, if $\psi$ as in
Lemma 2.3, then the dual wave solution is of the form
$\psi^+=e^{-kx-k^2y}\phi^+(Ux+Z,y,k)$, where
the coefficients $\xi_s^+(z+Z,y)$ of the formal series
\beq\label{psi2+}
\phi^+(z+Z,y,k)=e^{-by}\left(1+\sum_{s=1}^{\infty}\xi^+_s(z+Z,y)\, k^{-s}\right)
\eeq
are $\l$-periodic meromorphic functions of the variable $z\in \mathbb C^d$ with a
simple pole at the divisor $\Theta^{\,U}(y)$.

The ambiguity in the definition of $\psi$ does not affect the product
\beq\label{J0}
\psi^+\psi=\left(e^{-kx-k^2y}\Phi^{-1}\right)\left(\Phi e^{kx+k^2y}\right).
\eeq
Therefore, although each factor is only a local meromorphic function on
$\mathbb C^g\setminus \Sigma$, the coefficients $J_s$ of the product
\beq\label{J}
\psi^+\psi=\phi^+(Z,y,k)\phi(Z,y,k)=1+\sum_{s=2}^{\infty}J_s(Z,y)k^{-s}.
\eeq
are {\it global meromorphic functions} of $Z$. Moreover, the translational invariance
of $u$ implies that they have the form $J_s(Z,y)=J_s(Z+Vy)$. Each of the factors in
the left hand side of (\ref{J}) has a simple pole on $\Theta-Vy$. Hence, $J_s(Z)$ is
a meromorphic function on $X$ with a second order pole at $\Theta$.

From the definition of $\L$ it follows that
\beq\label{20}
\res_k\left(\psi^+(\L^n\psi)\right)=\res_k\left(\psi^+k^n\psi\right)=J_{n+1}.
\eeq
On the other hand, using the identity
\beq\label{dic}
\res_k \left(e^{-kx}\D_1\right)\left(\D_2e^{kx}\right)=\res_{\p}\left(\D_2\D_1\right),
\eeq
which holds for any two pseudo-differential operators (\cite{dickey}), we get
\beq\label{201}
\res_k(\psi^+\L^n\psi)=\res_k\left(e^{-kx}\Phi^{-1}\right)\left(\L^n\Phi e^{kx}\right)=
\res_{\p}\L^n=F_n.
\eeq
Therefore, $F_n=J_{n+1}$ and the lemma is proved.

Let $\bf {\hat F}$ be a linear space generated by $\{F_m, \ m=0,1,\ldots\}$, where we
set $F_0=1$. It is a subspace of the
$2^g$-dimensional space of the abelian functions that have at most second order pole at
$\Theta$. Therefore, for all but $\hat g=\dim\ {\bf \hat F}$ positive integers $n$,
there exist constants $c_{i,n}$ such that
\beq\label{f1}
F_n(Z)+\sum_{i=0}^{n-1} c_{i,n}F_i(Z)=0.
\eeq
Let $I$ denote the subset of integers $n$ for which there are no such constants. We call
this subset the gap sequence.
\begin{lem} Let $\L$ be the pseudo-differential operator corresponding to
a $\l$-periodic wave function $\psi$ constructed above. Then, for the differential operators
\beq\label{a2}
L_n=\L^n_++\sum_{i=0}^{n-1} c_{i,n}\L^{n-i}_+=0, \ n\notin I,
\eeq
the equations
\beq\label{lp}
L_n\,\psi=a_n(k)\,\psi, \ \ \ a_n(k)=k^n+\sum_{s=1}^{\infty}a_{s,n}k^{n-s}
\eeq
where $a_{s,n}$ are constants, hold.
\end{lem}
{\it Proof.} First note that from (\ref{lax}) it follows that
\beq\label{lax3}
[\p_y-\p_x^2+u,L_n]=0.
\eeq
Hence, if $\psi$ is a $\l$-periodic wave solution of (\ref{lax1})
corresponding to $Z\notin \Sigma$, then $L_n\psi$ is also a formal
solution of the same equation. That implies the equation
$L_n\psi=a_n(Z,k)\psi$, where $a$ is $\p_1$-invariant. The ambiguity in the definition of
$\psi$ does not affect $a_n$. Therefore, the coefficients of $a_n$ are well-defined
{\it global} meromorphic functions on $\mathbb C^g\setminus \Sigma$. The $\p_1$-
invariance of $a_n$ implies that $a_n$, as a function of $Z$, is holomorphic outside
of the locus. Hence it has an extension to a holomorphic function on $\mathbb C^g$.
Equations (\ref{tr1}) imply that $a_n$ is periodic with respect to the lattice
$\Lambda$. Hence $a_n$ is $Z$-independent. Note that $a_{s,n}=c_{s,n},\ s\leq n$.
The lemma is proved.

The operator $L_m$ can be regarded as a $Z\notin \Sigma$-parametric family   of
ordinary differential operators $L_m^Z$ whose coefficients have the form
\beq\label{lu}
L_m^Z=\p_x^n+\sum_{i=1}^m u_{i,m}(Ux+Z)\, \p_x^{m-i},\ \ m\notin I.
\eeq
\begin{cor} The operators $L_m^Z$
commute with each other,
\beq\label{com1}
[L_n^Z,L_m^Z]=0, \ Z\notin \Sigma.
\eeq
\end{cor}
From (\ref{lp}) it follows that $[L_n^Z,L_m^Z]\psi=0$. The commutator is an ordinary
differential operator. Hence, the last equation implies (\ref{com1}).

\begin{lem} Let $\A^Z,\ Z\notin \Sigma,$ be a commutative ring of ordinary differential
operators spanned by the operators $L_n^Z$. Then there is an irreducible algebraic
curve $\G$ of arithmetic genus $\hat g=\dim\ {\bf \hat F}$ such that $\A^Z$ is isomorphic
to the ring $A(\G,P_0)$ of the meromorphic functions on $\G$ with the only pole at
a smooth point $P_0$. The correspondence $Z\to \A^Z$ defines a holomorphic
imbedding of $X\setminus \Sigma$ into the space of torsion-free rank 1 sheaves $\F$ on $\G$
\beq\label{is}
j: X\backslash\Sigma\longmapsto \overline{\rm Pic}(\G).
\eeq
\end{lem}
{\it Proof.} It is the fundamental fact of the theory of commuting linear
ordinary differential operators (\cite{kr1,kr2,ch1,ch2,mum}) that there is a
natural correspondence
\beq\label{corr}
\A\longleftrightarrow \{\G,P_0, [k^{-1}]_1, \F\}
\eeq
between {\it regular} at $x=0$ commutative rings $\A$ of ordinary linear
differential operators containing a pair of monic operators of co-prime orders, and
sets of algebraic-geometrical data $\{\G,P_0, [k^{-1}]_1, \F\}$, where $\G$ is an
algebraic curve with a fixed
first jet $[k^{-1}]_1$ of a local coordinate $k^{-1}$ in the neighborhood of a smooth
point $P_0\in\G$ and $\F$ is a torsion-free rank 1 sheaf on $\G$ such that
\beq\label{sheaf}
H^0(\G,\F)=H^1(\G,\F)=0.
\eeq
The correspondence becomes one-to-one if the rings $\A$ are considered modulo conjugation
$\A'=g(x)\A g^{-1}(x)$.

Note that in \cite{kr1,kr2,ch1,ch2} the main attention was paid to the generic case of
the commutative rings corresponding to smooth algebraic curves.
The invariant formulation of the correspondence given above is due to Mumford \cite{mum}.

The algebraic curve $\G$ is called the spectral curve of $\A$.
The ring $\A$ is isomorphic to the ring $A(\G,P_0)$ of meromorphic functions
on $\G$ with the only pole at the puncture $P_0$. The isomorphism is defined by
the equation
\beq\label{z2}
L_a\psi_0=a\psi_0, \ \ L_a\in \A, \ a\in A(\G,P_0).
\eeq
Here $\psi_0$ is a common eigenfunction of the commuting operators. At $x=0$ it is
a section of the sheaf $\F\otimes\O(-P_0)$.

{\bf Important remark}. The construction of the correspondence (\ref{corr})
depends on a choice of initial point $x_0=0$. The spectral curve and the sheaf $\F$
are defined by the evaluations of the coefficients of generators of $\A$ and a finite
number of their derivatives at the initial point. In fact, the spectral curve
is independent on the choice of $x_0$, but the sheaf does depend on it, i.e. $\F=\F_{x_0}$.

Using the shift of the initial point it is easy to show that the correspondence
(\ref{corr}) extends to the commutative rings of operators whose coefficients are
{\it meromorphic} functions of $x$ at $x=0$. The rings of operators having poles at $x=0$
correspond to
sheaves for which the condition (\ref{sheaf}) is violated.

Let $\G^Z$ be the spectral curve corresponding to $\A^Z$.
Note, that due to the remark above,
it is well-defined for all $Z\notin \Sigma$.
The eigenvalues $a_n(k)$ of the operators $L_n^Z$ defined in (\ref{lp})
coincide with the Laurent expansions at $P_0$ of the meromorphic
functions $a_n\in A(\G^Z,P_0)$. They are $Z$-independent. Hence, the spectral curve
is $Z$-independent, as well, $\G=\G^Z$. The first statement of the lemma is thus
proven.

The construction of the correspondence
(\ref{corr}) implies that if the coefficients of operators $\A$
holomorphically depend on parameters then the algebraic-geometrical spectral data are
also holomorphic functions of the parameters. Hence $j$ is holomorphic away of $\Theta$.
Then using the shift of the initial point and the fact, that $\F_{x_0}$ holomorphically
depends on $x_0$, we get that $j$ holomorphically extends on $\Theta\setminus \Sigma$,
as well. The lemma is proved.

Recall, that a commutative ring $\A$ of linear ordinary differential operators
is called maximal if it is not contained in any bigger commutative ring.
Let us show that for a generic $Z$ the ring $\A^Z$ is maximal. Suppose that it is not.
Then there exits  $\a\in I$, where $I$ is the gap sequence defined above, such
that for each $Z\notin \Sigma$ there exists an operator $L_{\a}^Z$
of order $\a$ which commutes with $L_n^Z, n\notin I$. Therefore, it commutes with $\L$.
A differential operator commuting with $\L$ up to the order $O(1)$
can be represented in the
form $L_{\a}=\sum_{m<\a} c_{i,\a}(Z)\L^i_+$, where $c_{i,\a}(Z)$ are
$\p_1$-invariant functions of $Z$. It commutes with $\L$ if and only if
\beq\label{z3}
F_{\a}(Z)+\sum_{i=0}^{n-1} c_{i,\,\a}(Z)F_i(Z)=0,\ \ \p_1 c_{i,\,\a}=0.
\eeq
Note the difference between (\ref{f1}) and (\ref{z3}). In the first equation the
coefficients $c_{i,n}$ are constants. The $\l$-periodic wave solution of equation
(\ref{lax1}) is a common eigenfunction of all  commuting operators, i.e.
$L_{\a}\psi=a_{\a}(Z,k)\psi$, where
$a_{\a}=k^{\a}+\sum_{s=1}^{\infty} a_{s,\a}(Z)k^{\a-s}$ is $\p_1$-invariant.
The same arguments as those used in the proof of equation (\ref{lp}) show that
the eigenvalue $a_{\a}$ is $Z$-independent. We have
$a_{s,\,\a}=c_{s,\,\a},\ s\leq \a$. Therefore, the coefficients in (\ref{z3}) are
$Z$-independent. That contradicts the assumption that $\a\notin I$.

Our next goal is to prove finally the global existence of the wave function.
\begin{lem} Let the assumptions of the Theorem 1.2 hold. Then there exists a common
eigenfunction
of the corresponding commuting operators $L_n^Z$ of the form
$\psi=e^{kx}\phi(Ux+Z,k)$ such that
the coefficients of the formal series
\beq\label{psi6}
\phi(Z,k)=1+\sum_{s=1}^{\infty}\xi_s(Z)\, k^{-s}
\eeq
are global merormophic functions with a simple pole at $\Theta$.
\end{lem}
{\it Proof.} It is instructive to consider first the case when the spectral curve $\G$
of the rings $\A^Z$ is smooth. Then, as shown in (\cite{kr1, kr2}),
the corresponding common
eigenfunction of the commuting differential operators (the Baker-Akhiezer function),
normalized by the
condition $\psi_0|_{x=0}=1$, is of the form
(\cite{kr1, kr2})
\beq\label{ba}
\hat \psi_{0}={\hat\theta (\hat A(P)+\hat Ux+\hat Z)\,\hat\theta (\hat Z)\over
\hat\theta(\hat Ux+\hat Z)\,\hat\theta(\hat A(P)+\hat Z)}\,
e^{x\,\Omega(P)}.
\eeq
Here $\hat \theta (\hat Z)$ is the Riemann theta-function constructed with the help of the
matrix of $b$-periods of normalized holomorphic differentials on $\G$; $\hat A: \G\to J(\G)$
is the Abel map; $\Omega$ is the abelian integral corresponding to the second kind
meromorphic differential $d\Omega$ with the only pole of the form $dk$ at the puncture
$P_0$ and $2\pi i \hat U$ is the vector of its $b$-periods.

\noindent{\it Remark.}
Let us emphasize, that the formula (\ref{ba}) is not the result of solution of some
differential equations.
It is a direct corollary of analytic properties of the Baker-Akhiezer function
$\hat \psi_0(x,P)$ on the spectral curve:

$(i)$ {\it $\hat\psi_0$ is meromorphic function of $P\in \G\setminus P_0$; its pole
divisor
is of degree $\tilde g$ and is $x$-independent. It is non-special,
if the operators are regular at the normalization point $x=0$};

$(ii)$ {\it in the neighborhood of $P_0$ the function $\hat\psi_0$ has the form (\ref{ps})
(with $y=0$)}.

\noindent
From the  Riemann-Rokh theorem it follows that, if $\hat\psi_0$ exists, then it is unique.
It is easy to check that the function $\hat\psi_0$ given by (\ref{ba})
is single-valued on $\G$ and has all the desired properties.

The last factors in the numerator and the denominator of (\ref{ba}) are $x$-independent.
Therefore, the function
\beq\label{ba1}
\hat \psi_{BA}={\hat\theta (\hat A(P)+\hat Ux+\hat Z)\over
\hat\theta(\hat Ux+\hat Z)}\,
e^{x\,\Omega(P)}
\eeq
is also a common eigenfunction of the commuting operators.

In the neighborhood of $P_0$ the function $\hat \psi_{BA}$ has the form
\beq\label{sh10}
\hat \psi_{BA}=e^{kx}
\left(1+\sum_{s=1}^{\infty}{\tau_s (\hat Z+\hat Ux)\over \hat \theta(\hat U x+\hat Z)}\,k^{-s}
\right), \ \ k=\Omega,
\eeq
where $\tau_s(\hat Z)$ are global holomorphic functions.

According to Lemma 3.4, we have a holomorphic imbedding $\hat Z=j(Z)$
of $X\setminus\Sigma$ into $J(\G)$.
Consider the formal series $\psi=j^*\hat \psi_{BA}$.
It is globally well-defined out of $\Sigma$.
If $Z\notin \Theta$, then $j(Z)\notin \hat \Theta$
(which is the divisor on which the condition (\ref{sheaf}) is violated).
Hence, the coefficients of $\psi$ are regular out of $\Theta$. The singular locus
is at least of codimension 2. Hence, using once again Hartorg's arguments we can
extend $\psi$ on $X$.

If the spectral curve is singular, we can proceed along the same lines
using the generalization of (\ref{ba1})
given by the theory of Sato $\tau$-function
(\cite{wilson}). Namely, a set of algebraic-geometrical data (\ref{corr}) defines
the point of the Sato Grassmanian, and therefore, the corresponding $\tau$-function: $\tau(t;\F)$.
It is a holomorphic function of the variables $t=(t_1,t_2,\ldots)$,
and is a section of a holomorphic line bundle on $\overline{\rm Pic}(\G)$.

The variable $x$ is identified with the first time of the KP-hierarchy, $x=t_1$.
Therefore,
the formula for the Baker-Akhiezer function corresponding to a point of the Grassmanian
(\cite{wilson})
implies that the function $\hat \psi_{BA}$ given by the formula
\beq\label{baw}
\hat \psi_{BA}={\tau(x-k,-{1\over 2}k^2,-{1\over 3}k^3,\ldots; \F) \over
\tau(x,0,0,\ldots; \F)}e^{kx}
\eeq
is a common eigenfunction of the commuting operators defined by $\F$.
The rest of the arguments proving the lemma are
the same, as in the smooth case.
\begin{lem} The linear space $\bf {\hat F}$ generated by the abelian functions $\{F_0=1,
F_m=\res_{\p}\, \L^m\},$ is
a subspace of the space ${\bf H}$ generated by $F_0$ and by the abelian functions
$H_i=\p_1\p_{z_i}\ln \theta(Z)$.
\end{lem}
{\it Proof.} Recall that the functions $F_n$ are abelian
functions with at most second order pole on $\Theta$. Hence, a'priory
$\hat g=\dim\, {\bf \hat F}\leq 2^g.$ In order to prove the statement of the lemma
it is enough to show that $F_n=\p_1 Q_n$, where $Q_n$ is a
meromorphic function with a pole along $\Theta$. Indeed, if $Q_n$ exists,
then, for any vector $\l$ in the period lattice, we have $Q_n(Z+\l)=Q_n(Z)+c_{n,\l}$.
There is no abelian function with a simple pole on $\Theta$. Hence, there exists
a constant $q_n$ and two $g$-dimensional vectors $l_n,l_n'$, such that
$Q_n=q_n+(l_n,Z)+(l_n',h(Z))$, where $h(Z)$ is a vector with the coordinates
$h_i=\p_{z_i}\ln \theta$. Therefore, $F_n=(l_n,U)+(l_n',H(Z))$.

Let $\psi(x,Z,k)$ be the formal Baker-Akhiezer function defined in the previous lemma.
Then the coefficients $\varphi_s(Z)$ of the corresponding wave operator $\Phi$ (\ref{S})
are global meromorphic functions with poles on $\Theta$.

The left and  right action of pseudo-differential operators are formally adjoint,
i.e., for any two operators the equality $\left(e^{-kx}\D_1\right)\left(\D_2e^{kx}\right)=
e^{-kx}\left(\D_1\D_2e^{kx}\right)+\p_x\left(e^{-kx}\left(\D_3e^{kx}\right)\right)$
holds. Here $\D_3$ is a pseudo-differential operator whose coefficients are differential
polynomials in the coefficients of $\D_1$ and $\D_2$. Therefore, from (\ref{J0}-\ref{201})
it follows that
\beq\label{z8}
\psi^+\psi=1+\sum_{s=2}^{\infty}F_{s-1}k^{-s}=
1+\p_x\left(\sum_{s=2}^{\infty}Q_sk^{-s}\right).
\eeq
The coefficients of the series $Q$ are differential polynomials in the
coefficients $\varphi_s$ of the wave operator. Therefore, they are global meromorphic
functions of $Z$ with poles on $\Theta$. Lemma is proved.

In order to complete the proof of our main result we need one more standard fact
of the KP theory: flows of the KP hierarchy define deformations of the commutative
rings $\A$ of ordinary linear differential operators. The spectral curve
is invariant under these flows. For a given spectral curve $\G$ the orbits of
the KP hierarchy are isomorphic to the generalized Jacobian $J(\G)={\rm Pic}^0 (\G)$,
which is the equivalence classes of zero degree divisors on the spectral curve
(see details in \cite{shiota,kr1,kr2,wilson}).

The KP hierarchy in the Sato form is a system of  commuting differential
equation for a pseudo-differential operator $\L$
\beq\label{z4}
\p_{t_n}\L=[\L^n_+,\L]\,.
\eeq
If the operator $\L$ is as above. i.e., if it is defined by $\l$-periodic wave
solutions of equation (\ref{lax1}), then equations (\ref{z4}) are equivalent to
the equations
\beq\label{z5}
\p_{t_n}u=\p_x F_n.
\eeq
The first two times of the hierarchy are identified with the variables $t_1=x, t_2=y$.

Equations (\ref{z5}) identify the space $\hat {\bf F}_1$ generated by the functions
$\p_1 F_n$ with the tangent space of the KP orbit at  $\A^Z$.
Then, from Lemma 3.6 it follows that this tangent space is a subspace of the tangent space
of the abelian variety $X$. Hence, for any $Z\notin \Sigma$, the orbit of the KP flows
of the ring $\A^Z$ is in $X$, i.e. it defines an holomorphic imbedding:
\beq\label{imb}
i_Z:J(\G)\longmapsto X.
\eeq
From (\ref{imb}) it follows that $J(\G)$ is {\it compact}.

The generalized Jacobian of an algebraic curve is compact
if and only if the curve is {\it smooth} (\cite{mdl}). On a smooth algebraic curve
a torsion-free rank 1 sheaf is a line bundle, i.e. $\overline {\rm Pic} (\G)=J(\G)$.
Then (\ref{is}) implies that $i_Z$ is an isomorphism. Note, that
for the Jacobians of smooth algebraic curves the bad locus $\Sigma$ is empty
(\cite{shiota}), i.e. the imbedding $j$ in (\ref{is}) is defined everywhere
on $X$ and is inverse to $i_Z$. Theorem 1.2 is proved.

{\bf Acknowledgments}. The author would like to thank Enrico Arbarello for very
useful conversations and help during the preparation of this paper. The author is greatful
to Takahiro Shiota whose remarks helped the author to clarify some missing arguments in
the first version of the paper.

\end{document}